\documentclass{amsart}

\begin{document}

\title[Density Invariance of Operational Quantities]{Density Invariance of
Certain Operational Quantities of Bounded Linear Operators in Normed Spaces}

\author[M. K. Kinyon]{Michael K. Kinyon}
\thanks{The original version of this paper was written in the summer of
1995 while Ronald W. Cross was visiting the author at Indiana University
South Bend. The author wishes to thank Prof. Cross for sharing his insights
into linear operators and relations.}
\address{Department of Mathematics\\
Western Michigan University \\
Kalamazoo, MI  USA}
\email{mkinyon@wmich.edu}

\date{}

% theorem styles

\theoremstyle{plain}
\newtheorem{theorem}{Theorem}
\newtheorem{lemma}{Lemma}

% misc shortcut commands

\newcommand{\bR}{\mathbb{R}}
\newcommand{\cI}{\mathcal{I}}
\newcommand{\cP}{\mathcal{P}}
\newcommand{\cC}{\mathcal{C}}

\begin{abstract}
Based on ideas of R.W. Cross, a simplified proof is presented
of the density invariance of certain operational quantities associated with
bounded linear operators in normed vector spaces.
\end{abstract}

\maketitle

Let $X$, $Y$ denote normed linear spaces and let $T : X \to Y$ be a 
bounded linear operator.  Let $\cI(X)$ denote the collection of infinite 
dimensional subspaces of $X$.  For any subspace $M$ of $X$, let $S_M = 
\{ m \in M : \| m \| = 1\}$ and let $T|_M$ denote the restriction of $T$ to 
$M$.  We introduce the following quantities:

\begin{align*}
\Gamma (T)         &= \inf_{M \in \cI(X)} \| T|_M \| \\
\Delta (T)         &= \sup_{M \in \cI(X)} \Gamma (T|_M ) \\
\tau (T)           &= \sup_{M \in \cI(X)} \inf_{m \in S_M} \| Tm \| \\
\nabla (T)         &= \inf_{M \in \cI(X)} \tau (T|_M ) 
\end{align*}

These quantities were originally introduced in \cite{LS}, \cite{S}, and \cite{Se}
for bounded linear operators.  They were extended to unbounded linear operators
in \cite{Cr1} and \cite{Cr2}, and to linear relations in \cite{CP}.  See these
papers for applications of these quantities to the study of bounded and 
unbounded linear operators and relations.

A linear subspace $E$ of $X$ is said to be a {\em core} of $T$ if the graph 
of $T|_E$ is dense in the graph of $T$.  A quantity $f$ is said to be 
{\em densely invariant} if $f(T|_E ) = f(T)$ whenever $E$ is a core of $T$.

Our main result is the density invariance of the quantities $\Gamma$, $\Delta$, 
$\tau$, and $\nabla$.  For $\Gamma$, $\Delta$, and $\tau$, this result first 
appeared in \cite{Cr2}; for $\nabla$, the first appearance is in \cite{Cr3}.  
The purpose of this note is to give a simplified proof.  While the underlying
ideas are the same as those in \cite{Cr2} and \cite{Cr3}, the present proof 
has the following advantages: (1) the role played by the estimates is clearer, 
which considerably shortens the proof, and (2) the proof of the density 
invariance of $\Delta$ is direct and more obviously parallel to the proof of 
the density invariance of $\nabla$.

Our proof is for the case of bounded operators, but the case of unbounded
operators reduces to this case \cite{Cr2}, and the case of linear relations
reduces to the case of unbounded operators \cite{Cr3}.

We require the following \cite[Lemma IV.2.8(ii)]{Go}:

\begin{lemma}
Let $E$ be a dense linear subspace of $X$ and let $M$ be a closed,
finite codimensional subspace of $X$.  Then $E \cap M$ is dense in $M$.
\end{lemma}

\begin{theorem}
The quantities $\Gamma$, $\Delta$, $\tau$, and $\nabla$ are densely invariant.
\end{theorem}

\begin{proof}
Let $E$ be a core of $T$.  The case $\text{dim}\ X < \infty$ is trivial, so we 
assume that $X$ is infinite dimensional.  For $M \in \cI(X)$, we construct 
sequences $\{ m_n \} \subset M$ and $\{ x_n' \} \subset X'$ (here $X'$ is the 
dual of $X$) as follows.  Choose $m_1 \in S_M$ arbitrarily and choose $x_1'$ so 
that $\| x_1' \| = x_1' m_1 = 1$.  Suppose that $\{ m_1, \ldots, m_{n-1} \}$ and
$\{ x_1', \ldots, x_{n-1}' \}$ have been chosen.  Choose $m_n$ and $x_n'$ 
so that

\begin{equation}\label{eq1}
x_i' m_n = 0, \ \ \ 1 \leq i < n,
\end{equation}

\begin{equation}\label{eq2}
 \| m_n \| = \| x_n' \| = x_n' m_n = 1.
\end{equation}

Let $W = \text{span} \{ m_1, m_2, \ldots\}$.  It can be verified using 
(\ref{eq1}) and (\ref{eq2}) that the sequence $\{ m_n \}$ is linearly independent.
Thus $W \in \cI(X)$.  For an arbitrary element $w = \sum_1^n a_i m_i$, we shall 
verify that 

\begin{equation}\label{eq3}
 |a_i| \leq 2^{i-1} \| w \|, \ \ \ 1 \leq i \leq n.
\end{equation}

\noindent Now $a_1 = x_1'w$ by (\ref{eq1}), so $|a_1| \leq \| x_1' \| \| w \|
= \| w \|$.  Suppose that (\ref{eq3}) has been verified up to $i = j$.  From
(\ref{eq1}) and (\ref{eq2}), $x_{j+1}' w = \sum_1^j a_i x_{j+1}' m_i + a_{j+1}$.
By the induction hypothesis, 
\[
|a_{j+1}| \leq | x_{j+1}' w | + \sum_1^j |a_i| |x_{j+1}' m_i | \leq \| w \| +
\sum_1^j 2^{i-1} \| w \| \leq 2^j \| w \|
\]
as required.

Let $\epsilon$ be given with $0 < \epsilon < 1$.  For each $n$, the intersection
of the subspace $E$ with the intersection $N$ of the kernels of the bounded linear 
functionals $x_1', \ldots, x_n'$ is dense in $N$ by the lemma, and hence there
exists a sequence $\{z_n\}$ in $E$ such that

\begin{equation}\label{eq4}
 \| z_n - m_n \| \leq 2^{1-2n} \epsilon \ \text{min} \{1, c/\|T\| \}
\end{equation}

\noindent and

\begin{equation}\label{eq5}
 x_i'z_n = 0, \ \ \ 1 \leq i < n
\end{equation}

\noindent where $c > 0$ is an arbitrary constant.

Let $L = \text{span} \{ z_1, z_2, \ldots\}$.  The linear independence of the
sequence $\{ z_n \}$ follows readily from (\ref{eq4}), (\ref{eq5}), and the
linear independence of $\{ m_n \}$.  Thus $L \in \cI(E)$.

Define a linear bijection $A : L \to W$ by $A z_i = m_i$.  For $z = \sum_1^n 
a_i z_i \in L$, we have

\begin{equation}\label{eq6}
   \begin{split}
    \| z - Az \| &\leq \sum_1^n |a_i| \| z_i - A z_i \| \leq \sum_1^n 2^{i-1} 
    \| Az \| \| z_i - A z_i \| \\
                 &\leq \epsilon \sum_1^n 2^{i-1} 2^{1-2i}\ \text{min}\{1,c/\|T\|\} 
    \|Az\| < \epsilon \ \text{min}\{1,c/\|T\|\} \| Az \|
   \end{split}
\end{equation}   

\noindent by (\ref{eq3}) and (\ref{eq4}).  Consequently

\begin{equation}\label{eq7}
   \begin{split}
    (1 - \epsilon) \| Az \| &\leq (1 - \epsilon \ \text{min}\{1,c/\|T\|\}) \| Az\| \\
                            &< \| Az \| - \| z - Az \| \\
                            &\leq \| z \| \\
                            &\leq \| Az \| + \| z - Az \| \\
                            &< (1 + \epsilon \ \text{min}\{1,c/\|T\|\}) \| Az\| \\
                            &\leq ( 1 + \epsilon) \| Az \|.
   \end{split}
\end{equation}   
    
For $z \in L$, we have

\begin{align*}
    \| Tz \| &\geq \| TAz \| - \| T ( z - Az )\| \\
             &\geq \| TAz \| - \| T \| \| z - Az \| \\
             &> \| TAz \| - \|T\| \epsilon \ \text{min}\{1,c/\|T\|\} \|Az\| \\
             &\geq \| TAz \| - \epsilon c \| Az \|.
\end{align*}             

\noindent by (\ref{eq6}).  Thus using (\ref{eq7}),

\begin{equation}\label{eq8}
\frac{\| Tz \|}{\| z\|} > ( \frac{\| TAz \|}{\| Az\|} - \epsilon c)
(1 + \epsilon )^{-1}.
\end{equation}

\noindent Similarly, we have for $z \in L$,

\begin{align*}
    \| Tz \| &\leq \| TAz \| + \| T ( z - Az )\| \\
             &\leq \| TAz \| + \| T \| \| z - Az \| \\
             &< \| TAz \| + \|T\| \epsilon \ \text{min}\{1,c/\|T\|\} \|Az\| \\
             &\leq \| TAz \| + \epsilon c \| Az \|
\end{align*}             

\noindent by (\ref{eq6}).  Thus using (\ref{eq7}),

\begin{equation}\label{eq9}
\frac{\| Tz \|}{\| z\|} < ( \frac{\| TAz \|}{\| Az\|} + \epsilon c)
(1 - \epsilon )^{-1}.
\end{equation}

\medskip

(a) The quantity $\Gamma$: We clearly have $\Gamma (T) \leq \Gamma (T|_E )$.
Suppose that $\Gamma (T) < \Gamma (T|_E ) - \delta$ where $\delta > 0$.  Set
$c = \Gamma (T|_E ) - \delta$.  Choose $M \in \cI(X)$ such that

\begin{equation}\label{eq10}
\| T|_M \| < c.
\end{equation}

\noindent Then from (\ref{eq9}) and (\ref{eq10})
\[
\frac{\| Tz \|}{\| z\|} < \frac{1 + \epsilon}{1 - \epsilon} c
\]

\noindent for $z \in L$.  Hence
\[
\Gamma ( T|_E ) = \inf_{V \in \cI(E)} \| T|_V \| \leq \| T|_L \| \leq 
\frac{1 + \epsilon}{1 - \epsilon} c.
\]
Since $\epsilon$ is arbitrary, $\Gamma ( T|_E ) \leq c$, a contradiction.

\medskip

(b) The quantity $\tau$: We clearly have $\tau (T) \geq \tau (T|_E )$.
Suppose that $\tau (T) > \tau (T|_E ) + \delta$ where $\delta > 0$.  Set
$c = \tau (T|_E ) + \delta$.  Choose $M \in \cI(X)$ such that

\begin{equation}\label{eq11}
\frac{\|Tm\|}{\|m\|} > c \ \text{for}\ m \in M.
\end{equation}

\noindent Then from (\ref{eq8}) and (\ref{eq11})
\[
\frac{\| Tz \|}{\| z\|} > \frac{1 - \epsilon}{1 + \epsilon} c
\]

\noindent for $z \in L$.  Hence
\[
\tau ( T|_E ) \geq \inf_{z \in L} \frac{\| Tz \|}{\| z\|} \geq 
\frac{1 - \epsilon}{1 + \epsilon} c.
\]
Since $\epsilon$ is arbitrary, $\tau ( T|_E ) \geq c$, a contradiction.

\medskip

(c) The quantity $\Delta$: We clearly have $\Delta (T|_E ) \leq \Delta (T)$.
Let $\delta > 0$ be given and choose $M_0 \in \cI(X)$ such that
$\Gamma (T|_{M_0}) > \Delta (T) - \delta$.  Set $c = \Delta (T) - \delta$.  
Then if $M \in \cI(M_0)$,

\begin{equation}\label{eq12}
\| T|_M \| > c.
\end{equation}

From (\ref{eq8}) we have
\[
\| T|_V \| \geq ( \| T|_{AV} \| - \epsilon c)(1 + \epsilon )^{-1}
\]
for $V \in \cI(L)$.  Since $AV \in \cI(M) \subset \cI(M_0)$, it follows from
(\ref{eq12}) that 
\[
\| T|_V \| \geq \frac{1 - \epsilon}{1 + \epsilon} c
\]
for $V \in \cI(L)$.  Thus
\[
\Gamma (T|_L ) = \inf_{V \in \cI(L)} \| T|_V \| \geq 
\frac{1 - \epsilon}{1 + \epsilon} c.
\]
Taking the supremum of the left side over $L \in \cI(E)$, we have
\[
\Delta (T|_E ) \geq \frac{1 - \epsilon}{1 + \epsilon} c = 
\frac{1 - \epsilon}{1 + \epsilon} ( \Delta (T) - \delta ).
\]
Since $\epsilon$ and $\delta$ are arbitrary, it follows that 
$\Delta (T|_E ) \geq \Delta (T)$.

\medskip

(d) The quantity $\nabla$: We clearly have $\nabla (T|_E ) \geq \nabla (T)$.
Let $\delta > 0$ be given and choose $M_0 \in \cI(X)$ such that
$\tau (T|_{M_0}) < \nabla (T) + \delta$.  Set $c = \nabla (T) + \delta$.  
Then if $M \in \cI(M_0)$,

\begin{equation}\label{eq13}
\inf_{m \in S_M} \| Tm \| < c.
\end{equation}

From (\ref{eq9}) we have
\[
\inf_{z \in V}\frac{\| Tz \|}{\|z\|} \leq ( \inf_{Az \in AV}
\frac{\| TAz \|}{\| Az\|} + \epsilon c)(1 - \epsilon )^{-1}
\]
for $V \in \cI(L)$.  Since $AV \in \cI(M) \subset \cI(M_0)$, it follows from
(\ref{eq13}) that 
\[
\inf_{z \in V}\frac{\| Tz \|}{\|z\|} \leq \frac{1 + \epsilon}{1 - \epsilon} c
\]
for $V \in \cI(L)$.  Thus
\[
\tau (T|_L ) = \sup_{V \in \cI(L)} \inf_{z \in V}\frac{\| Tz \|}{\|z\|} \leq 
\frac{1 + \epsilon}{1 - \epsilon} c.
\]
Taking the infimum of the left side over $L \in \cI(E)$, we have
\[
\nabla (T|_E ) \leq \frac{1 + \epsilon}{1 - \epsilon} c = 
\frac{1 + \epsilon}{1 - \epsilon} ( \nabla (T) + \delta ).
\]
Since $\epsilon$ and $\delta$ are arbitrary, it follows that 
$\nabla (T|_E ) \leq \nabla (T)$.
\end{proof}

\end{document}